\documentclass[12pt,reqno]{amsart}

\usepackage{a4wide}
\usepackage{amsmath}
\usepackage{amsthm}
\usepackage{amsfonts}
\usepackage{amssymb}
\usepackage{graphicx}
\usepackage[dvipsnames]{xcolor}
\usepackage{cite}
\usepackage{mathscinet}

\usepackage{color,soul}
\usepackage{xcolor}

\usepackage{todonotes}

\usepackage[bookmarksnumbered, colorlinks, plainpages]{hyperref}

\theoremstyle{plain}
\newtheorem*{thma}{Theorem A}
\newtheorem{thm}{Theorem}
\newtheorem{lem}[thm]{Lemma}

\newtheorem{prop}[thm]{Proposition}

\theoremstyle{definition}

\newtheorem{defi}[thm]{Definition}

\def\supp{\operatorname{supp}}
\def\esup{\operatornamewithlimits{ess\,sup}}

\def\RHS{\operatorname{RHS}}
\def\LHS{\operatorname{LHS}}

\newcommand{\M}{\mathcal{M}}
\newcommand{\Z}{\mathbb{Z}}
\newcommand{\Mpl}{\M^+}

\allowdisplaybreaks

\begin{document}

\author[Amiran Gogatishvili et al.]{Amiran Gogatishvili, Zden\v ek Mihula, Lubo\v s Pick, Hana Tur\v cinov\'a \and Tu\u{g}\c{c}e \"{U}nver}
	
\title[Weighted inequalities for a superposition of operators]{Weighted inequalities for a superposition of the Copson operator and the Hardy operator}

\address{Amiran Gogatishvili,
 Institute of Mathematics of the
 Czech Academy of Sciences,
 \v Zitn\'a~25,
 115~67 Praha~1,
 Czech Republic, ORCID 0000-0003-3459-0355}
\email{gogatish@math.cas.cz}

\address{Zden\v ek Mihula, Czech Technical University in Prague, Faculty of Electrical Engineering, Department of Mathematics, Technick\'a~2, 166~27 Praha~6, Czech Republic; AND Department of Mathematical Analysis, Faculty of Mathematics and Physics, Charles University, Sokolovsk\'a~83, 186~75 Praha~8,	Czech Republic, ORCID 0000-0001-6962-7635}
\email{mihulzde@fel.cvut.cz \& mihulaz@karlin.mff.cuni.cz}

\address{Lubo\v s Pick,
 Department of Mathematical Analysis,
	Faculty of Mathematics and Physics,
	Charles University,
	Sokolovsk\'a~83,
	186~75 Praha~8,
	Czech Republic, ORCID 0000-0002-3584-1454}
\email{pick@karlin.mff.cuni.cz}

\address{Hana Tur\v cinov\'a,
 Department of Mathematical Analysis,
	Faculty of Mathematics and Physics,
	Charles University,
	Sokolovsk\'a~83,
	186~75 Praha~8,
	Czech Republic, ORCID 0000-0002-5424-9413}
\email{turcinova@karlin.mff.cuni.cz}

\address{T. \"{U}nver (\textit{corresponding author}),
Institute of Mathematics of the
Czech Academy of Sciences,
\v Zitn\'a~25,
115~67 Praha~1,
Czech Republic \newline
Faculty of Science and Arts,
Kirikkale University,
71450 Yahsihan, Kirikkale,
Turkey, ORCID 0000-0003-0414-8400}
\email{unver@math.cas.cz \& tugceunver@kku.edu.tr}
	

\thanks{This research was supported in part by the grant P201-18-00580S of the Czech Science Foundation and by the Danube Region Grant no.~8X2043 of the Czech Ministry of Education, Youth and Sports. The research of  A.~Gogatishvili was also supported by Czech Academy of Sciences RVO: 67985840. The research of Z.~Mihula was supported by the project OPVVV CAAS CZ.02.1.01/0.0/0.0/16\_019/0000778 and by the grant SVV-2020-260583. The research of L. Pick was supported in part by the grant P201/21-01976S of the Czech Science Foundation. The research of H.~Tur\v cinov\'{a} was supported in part by the Grant Schemes at Charles University, reg. no. CZ.02.2.69/0.0/0.0/19\_073/0016935, by the Primus research programme PRIMUS/21/SCI/002 of Charles University, by the grant SVV-2020-260583 and Charles University Research program No.~UNCE/SCI/023. The research of T.~\"{U}nver was supported by the grant of The Scientific and Technological Research Council of Turkey (TUBITAK), Grant No: 1059B192000075. Part of the work on this project was carried out during the meeting Per Partes held at Horn\'{\i} Lyse\v ciny, June 2--6, 2021.}

\begin{abstract}
We study a three-weight inequality for the superposition of the Hardy operator and the Copson operator, namely
\begin{equation*}
	\bigg(\int_a^b \bigg(\int_t^b \bigg(\int_a^s f(\tau)^p v(\tau) \,d\tau \bigg)^\frac{q}{p} u(s) \,ds \bigg)^{\frac{r}{q}} w(t) \,dt  \bigg)^{\frac{1}{r}} \leq C \int_a^b f(t)\,dt,
\end{equation*}
in which $(a,b)$ is any nontrivial interval, $q,r$ are positive real parameters and $p\in(0,1]$. A simple change of variables can be used to obtain any weighted $L^p$-norm with $p\ge1$ on the right-hand side. Another simple change of variables can be used to equivalently turn this inequality into the one in which the Hardy and Copson operators swap their positions. We focus on characterizing those triples of weight functions $(u,v,w)$ for which this inequality holds for all nonnegative measurable functions $f$ with a constant independent of $f$.

We use a new type of approach based on an innovative method of discretization which enables us to avoid duality techniques and therefore to remove various restrictions that appear in earlier work.

This paper is dedicated to Professor Stefan Samko on the occasion of his 80th birthday.
\end{abstract}
	
\bibliographystyle{abbrv}
	
\maketitle

\section{Introduction and the main result}

The main purpose of this paper is to introduce a new line of argument which enables one to obtain a previously unavailable characterization of the validity of certain specific inequalities involving superposition of integral operators of Copson and Hardy type and three weight functions.

More precisely, given $a,b\in[-\infty,\infty]$, $a<b$, and parameters $q,r \in(0,\infty)$ and $p\in(0,1]$, we characterize all triples $(u,v,w)$ of \textit{weights} (i.e.~positive measurable functions) on $(a,b)$ such that there exists a constant $C>0$ with which the inequality
\begin{equation} \label{main-iterated-intro}
	\bigg(\int_a^b \bigg(\int_t^b \bigg(\int_a^s f(\tau)^p v(\tau) \,d\tau \bigg)^\frac{q}{p} u(s) \,ds \bigg)^{\frac{r}{q}} w(t) \,dt  \bigg)^{\frac{1}{r}} \leq C \int_a^b f(t) \,dt
\end{equation}
holds for every nonnegative measurable function $f$ on $(a,b)$. Let us note that the restriction $p\in(0,1]$ is natural and does not cause any weakness. Indeed, the inequality is obviously impossible without it as, if $p>1$, one can always easily construct a function $f$ that makes the integral on the left diverge while keeping the right-hand side finite.

The inequality~\eqref{main-iterated-intro}, as a certain ``mother figure'', immediately paves the way to many other important inequalities. For instance, one can easily swap the order of the two inner integral operators and obtain the inequality
\begin{equation} \label{main-iterated-intro-swap}
	\bigg(\int_a^b \bigg(\int_a^t \bigg(\int_s^b f(\tau)^p v(\tau)\,d\tau \bigg)^\frac{q}{p} u(s) \,ds \bigg)^{\frac{r}{q}} w(t) \,dt  \bigg)^{\frac{1}{r}} \leq C \int_a^b f(t) \,dt
\end{equation}
instead of~\eqref{main-iterated-intro}. This is achieved by a simple change of variables $\tau\mapsto -\tau$ in the innermost integral on the left side of~\eqref{main-iterated-intro} and following the forced changes from thereafter. Similarly, one can turn~\eqref{main-iterated-intro} to the inequality
\begin{equation} \label{main-iterated-intro-rhs}
	\bigg(\int_a^b \bigg( \int_t^b \bigg( \int_a^s f(\tau)\,d\tau \bigg)^q u(s) \,ds \bigg)^{\frac{r}{q}} w(t) \,dt  \bigg)^{\frac{1}{r}} \leq C \left(\int_a^b f(t)^pv(t) \,dt\right)^{\frac{1}{p}}
\end{equation}
with $p\ge1$ by performing the replacements (in this order) $f\mapsto f^{\frac{1}{p}}v^{-\frac{1}{p}}$, $v\mapsto v^{-p}$, $q\mapsto qp$, $r\mapsto rp$ and, finally, $p\mapsto \frac{1}{p}$, in~\eqref{main-iterated-intro}. Our characterization of~\eqref{main-iterated-intro} thus immediately yields necessary and sufficient conditions for ~\eqref{main-iterated-intro-swap},~\eqref{main-iterated-intro-rhs}, and their various combinations.

The key innovation is contained in methods of proofs which are based on new discretization techniques that require neither duality methods nor nondegeneracy conditions on weights.

In the theory of weighted inequalities, questions involving iterations of operators have recently been constituting the cutting edge. The subject has been rather fashionable for some time, mainly because inequalities involving compositions of operators, on the one hand, have an impressive array of important applications, while, on the other hand, are quite difficult to handle.

There is plenty of motivation for studying weighted inequalities for a composition of operators, and it pours down from various sources, rather different in spirit. A notable one is the theory of Sobolev-type embeddings where, during the last two decades, various forms of the so-called reduction principles have flourished. The reduction principle is a powerful method which establishes an, perhaps somewhat surprising, equivalence between a difficult problem involving differential operators in several variables, such as a Sobolev-type embedding, and a weighted inequality for an integral operator acting on functions defined on an interval. For the first-order embedding, this is usually achieved by an effective use of some sort of the P\'olya--Szeg\H o principle, and the resulting operator is then always a weighted Copson operator.

For higher-order embeddings, however, the P\'olya--Szeg\H o principle does not work, and one needs some new way of argumentation. Here, once again, approaches vary. For Euclidean--Sobolev embeddings (in which functions are defined on a sufficiently regular subdomain of the Euclidean ambient space $\mathbb R^n$ endowed with the Lebesgue measure), an effective use of interpolation theory leads to satisfactory results \cite{Ke-Pi:06}. However, when the underlying domain is not sufficiently regular, or, for instance, when $\mathbb R^n$ is replaced by the Gauss space $(\mathbb R^n,\gamma_n)$, which is still $\mathbb R^n$, but endowed with the Gauss probability measure
\begin{equation*}
    \gamma_n(x) = (2\pi)^{-\frac{n}{2}}e^{-\frac{|x|^2}{2}}\,dx,
\end{equation*}
interpolation methods turn out to be ineffective and have to be replaced by something else. An extremely efficient argument in this matter, based on the isoperimetric inequality combined with an iteration technique, was developed in \cite{Ci-Pi-Sl:15}, and later applied again, this time to higher-order trace embeddings, in \cite{Ci-Pi:16}.

Now here comes the interesting part. While for some specific weights (typically for power weights that occur in reduction principles for Euclidean--Sobolev embeddings) the iterated operator is pointwise equivalent to a~suitable single Copson operator \cite{Ke-Pi:06}, this is impossible in general (for example for Euclidean--Sobolev embeddings on bad enough domains or for Gaussian--Sobolev embeddings this approach fails \cite{Ci-Pi:09}). In result, after performing the reduction principle, one has to grapple either with a~kernel operator, or with a~superposition of two or more operators \cite{Ci-Pi-Sl:15}. Specifically, superposition of integral operators of different type (Hardy vs.~Copson) is encountered for example when operators obtained from the reduction principle are applied to one of many operator-induced function spaces. The simplest examples of these are spaces whose norm involves the operation of the maximal nonincreasing rearrangement such as weak spaces, Marcinkiewicz spaces, etc., but there are more sophisticated ones which are also important.

Another variety of applications in a completely different direction can be found in the theory of function spaces and interpolation theory. These shelter, among others, questions concerning sharp embeddings between important structures \cite{GAG,Pe:21}, K\"othe duals of function spaces \cite{GKPS,Tu:21,Go-Jo-Ok-Pe:07}, inequalities restricted to cones of functions such as those of monotone or concave functions \cite{GOP,Gr:98,GogStep}, or inequalities involving bilinear and multilinear operators \cite{Be-Or:17}. The blocking technique appearing in~\cite{Gr:98} was also independently developed and applied for integral Hardy-type inequalities in~\cite{Go:98}.

Several results were obtained recently for iterations of operators of identical type, however always under some rather unpleasant restrictions.

One of the earliest treatments of iteration of identical operators was most likely carried out in~\cite{Bu:13}. The authors consider an $n$-dimensional problem
and using radial weights they reduce it to the iterated Hardy inequality and treat some particular cases of parameters by discretization. Later in~\cite{Go:17},
inequalities involving Hardy-Hardy or Hardy-Copson iteration were fully characterized. Owing to the fact that a reduction technique was used, the characterization obtained is more complicated and nonstandard. In~\cite{Ka:19} both cases of iterations (Hardy-Hardy and Hardy-Copson) involving a kernel and using a different discretization are
considered, in the simplest case of parameters, and characterization is obtained. Recently, iteration of Copson operators was treated in~\cite{Kr-Pi:20}, restricted to nondegenerate weights.

Inequalities for superposition of the Copson and Hardy operators were studied in~\cite{GMP}; however, the
results obtained there were restricted to nondegenerate weights. Next, in~\cite{Oi-Ka:08}, a three-weight inequality was characterized, motivated by a specific inequality in which a weighted norm of
a mean value is compared to that of the derivative of a given function. Techniques of proofs in that paper are related to~\cite{Si-St:96}. The result was later
revisited several times, see e.g.~\cite{Be-Or:17,Oi-Ka:15}, where also further
applications to bilinear operators are pointed out. Particular cases and related topics had been studied earlier, see for
instance~\cite{Ev-Go-Op:18} or \cite{Go:02} for $p=1$, or \cite{GOP} for $p=\infty$, or~\cite{Si-St:96} and~\cite{Go-Jo-Ok-Pe:07} for special cases of weights.
The subject is also intimately related to the new type of spaces governed by operator-induced norms that have been appearing recently in connection with various other tasks, notably from embeddings of Sobolev spaces endowed with slowly-decaying upper Ahlfors measures \cite{Ci-Pi-Sl:20,Ci-Pi-Sl:20a,Tu:21}.

Let us recall that discretization techniques have been around for some time. In the late 1980's and early 1990's they proved to be very useful for example in the theory of one-sided operators and ergodic theory, see e.g.~\cite{Sa:86,Ma:89,Sa:90,Ma:93} and all the huge amount of subsequent work. In the early 2000's, they were used in order to solve some problems in the theory of classical Lorentz spaces that had been open for long time, see~\cite{GoPi:03,CGMP2}. Later various authors spent considerable efforts in order to chip away certain technical obstacles such as nondegeneracy assumptions with varying success - consider e.g.~\cite{Ev-Go-Op:18} or~\cite{KMT} and the references therein. However, this research is far from being complete.

Let us note that the current paper is closely related to the project \cite{Go-Pi-Un:21} currently under preparation, in which some of the new discretization methods presented here will be applied to a different problem, namely to an inequality involving the Hardy operator on one side and the Copson operator on the other, cf.~\cite{CGMP2}.

We shall now present our principal result, that is, a complete characterization of~\eqref{main-iterated-intro}. We will formulate it in the form of a single theorem. We shall need the following notation. For $a,b\in[-\infty,\infty]$, $a<b$, and $p\in(0,1]$, let
\begin{align*}
    V_p(a, b) :=
    \begin{cases}
        \big(\int_a^b v^{\frac{1}{1-p}}\big)^{\frac{1-p}{p}} & \text{if $0<p<1$,}
            \\
        \esup\limits_{t \in (a, b)} v(t) & \text{if $p=1$.}
    \end{cases}
\end{align*}

Our main result is:

\begin{thma}
Let $a,b\in[-\infty,\infty]$, $a<b$, $q,r \in(0,\infty)$, $p\in(0,1]$, and let $u, v, w$ be weights on $(a,b)$. Then there exists a constant $C>0$ such that the inequality \eqref{main-iterated-intro} holds for all nonnegative measurable functions $f$ on $(a, b)$ if and only if one of the following conditions is satisfied:
	
\rm(i) $1 \le r $, $1\leq q$,
\begin{equation*}
C_1 :=  \sup_{ t \in (a, b)} \bigg(\int_a^t w(s) \,ds\bigg)^{\frac{1}{r}} \esup_{s \in (t,b)} \bigg(\int_s^b u(\tau) \,d\tau\bigg)^{\frac{1}{q}} V_p(a, s) < \infty
\end{equation*}
and
\begin{equation*}
C_2 := \sup_{ t \in (a, b)} \bigg(\int_t^b w(s) \bigg(\int_s^b u(\tau) \,d\tau \bigg)^{\frac{r}{q}} \,ds \bigg)^{\frac{1}{r}} V_p(a, t) < \infty;
\end{equation*}

\rm(ii) $ 1 \le r $, $q<1$, $C_2 < \infty$ and
\begin{equation*}
C_3 :=\sup_{t\in (a, b)} \bigg(\int_a^t w(s) \,ds \bigg)^{\frac{1}{r}}  \bigg(\int_{t}^{b} \bigg(\int_s^b u(\tau) \,d\tau \bigg)^{\frac{q}{1-q}} u(s) V_p(a, s)^{\frac{q}{1-q}} \,ds \bigg)^{\frac{1-q}{q}} < \infty;
\end{equation*}

\rm(iii) $r < 1$, $1 \leq q$,
\begin{equation*}
C_4 :=\bigg(\int_a^b  \bigg(\int_a^t w(s) \,ds \bigg)^{\frac{r}{1-r}} w(t) \esup_{s \in (t, b)} \bigg(\int_s^{b} u(\tau) \,d\tau \bigg)^{\frac{r}{q(1-r)}} V_p(a, s)^{\frac{r}{1-r}} \,dt \bigg)^{\frac{1-r}{r}} <\infty
\end{equation*}
and
\begin{align*}
C_5 & := \bigg(\int_a^b  \bigg(\int_t^b w(s) \bigg(\int_s^b u(\tau) \,d\tau \bigg)^{\frac{r}{q}} \,ds \bigg)^{\frac{r}{1-r}} w(t)
    \bigg(\int_t^{b} u(\tau) \,d\tau \bigg)^{\frac{r}{q}} V_p(a, t)^{\frac{r}{1-r}} \,dt \bigg)^{\frac{1-r}{r}} < \infty;
\end{align*}

\rm(iv) $r < 1$, $q < 1$, $C_5 < \infty$ and
\begin{equation*}
C_6 :=\bigg( \int_a^b \bigg(\int_a^t w(s) \,ds\bigg)^{\frac{r}{1-r}} w(t) \bigg( \int_{t}^{b} \bigg(\int_{s}^{b} u(\tau) \,d\tau \bigg)^{\frac{q}{1-q}} u(s) V_p(a,s)^{\frac{q}{1-q}}\,ds \bigg)^{\frac{r(1-q)}{q(1-r)}} \,dt \bigg)^{\frac{1-r}{r}} < \infty.
\end{equation*}
Moreover, the best constant $C$ in the inequality \eqref{main-iterated-intro} satisfies
\begin{equation}\label{thm:optimal_constant}
    C\approx
    \begin{cases}
        C_1 + C_2 &\text{in the case (i),}
            \\
        C_2 + C_3 &\text{in the case (ii),}
            \\
        C_4 + C_5 &\text{in the case (iii),}
            \\
        C_5 + C_6 &\text{in the case (iv).}
    \end{cases}
\end{equation}
\end{thma}

The proof is based on a new type of discretization which avoids the use of any kind of duality principle, enabling us thereby to obtain the result in the required generality.

Theorem~A is proved in Section~\ref{S:proofs}, along with a side theorem which gives another characterization of~\eqref{main-iterated-intro}. Key ingredients of the proofs are collected in Section~\ref{S:discrete}.

\section{Background discretization results}\label{S:discrete}

In this section we shall establish the background discretization material that will be needed in the proof of the main result. We first fix notation and conventions used in this paper. We denote by $\LHS(*)$ and $\RHS(*)$ the left-hand side and right-hand side of the inequality numbered by $*$, respectively. We adhere to the usual convention that $\frac1{\infty}=0\cdot\infty=\frac{\infty}{\infty}=\frac{0}{0}=0$. We denote by $\Mpl(c,d)$ the set of all nonnegative measurable functions on $(c,d)$. By increasing we mean strictly increasing. Finally, the small letters $i$ and $k$ are always integers, which are reserved for indices. In particular, when we write $N\leq k\leq M$, in which $N$ and $M$ can be $-\infty$ and $\infty$, respectively, we mean $k\in\Z$, $N\leq k\leq M$. This convention is accordingly modified for similar inequalities and the index $i$ in the obvious way.
\begin{defi}
Let $N \in \Z\cup\{-\infty\}$, $M \in \Z\cup\{+\infty\}$, $N<M$, and $\{a_k\}_{k=N}^M$ be a sequence of positive numbers. We say that $\{a_k\}_{k=N}^M$ is \textit{strongly increasing} if
\begin{equation*}
    \inf\bigg\{ \frac{a_{k+1}}{a_k},\quad N\leq k < M\bigg\} > 1.
\end{equation*}
\end{defi}

Our approach is based on a fine discretization of the inequality in question. Before we start doing that, we need some new information of a general kind from the discrete world. The following lemma is contained in the manuscript \cite{Go-Pi-Un:21}, where it is also proved. However, since the manuscript is not publicly available yet, we include its proof here for the reader's convenience.

\begin{lem}
Let $s>0$, $M \in \Z\cup \{+\infty\}$. Assume that $\{a_k\}_{k=-\infty}^M$ and $\{b_k\}_{k=-\infty}^M$ are sequences of nonnegative numbers such that $\{b_k\}_{k=-\infty}^M$ is nondecreasing. Then
\begin{equation}\label{difference-u}
    \sum_{k=-\infty}^{M} a_k \bigg(\sum_{i=k}^{M} a_i \bigg)^s b_k \approx \sum_{k=-\infty}^{M} (b_k  - b_{k-1} ) \bigg(\sum_{i= k}^{M} a_i\bigg)^{s+1} + \bigg(\sum_{k=-\infty}^{M} a_k\bigg)^{s+1} \lim_{k\rightarrow -\infty} b_k,
\end{equation}
in which the multiplicative constants depend only on $s$.
\end{lem}

\begin{proof}
First, assume that $\lim_{k\to-\infty}b_k > 0$. Thanks to this assumption, we have that
\begin{equation*}
\lim_{N\to -\infty}\bigg(\sum_{k=N}^M a_k\bigg)^{s+1} b_N = \bigg(\sum_{k=-\infty}^{M} a_k\bigg)^{s+1} \lim_{k\rightarrow -\infty} b_k,
\end{equation*}
whether the series converges or diverges. Let $N\in\Z$, $N<M$. By virtue of Abel's lemma, we have that
 \begin{equation}\label{L:equiv.cond.:Abel}
     \sum_{k=N}^M c_k b_k =  \sum_{k=N+1}^M (b_k  - b_{k-1}) \sum_{i= k}^M c_i + \bigg(\sum_{k=N}^M c_k\bigg) b_N
 \end{equation}
for every sequence $\{c_k\}_{k=N}^M$ of nonnegative numbers. Set $c_k = a_k \big(\sum_{i=k}^M a_i\big)^s$ for $k\in\Z$, $N\leq k\leq M$. Applying power rules (cf.~e.g.~\cite[Lemma~1 and Lemma~1']{Be-Gr:05}), we get
 \begin{align*}
 \sum_{k=N}^M a_k \bigg(\sum_{i=k}^M a_i \bigg)^s b_k & =  \sum_{k=N+1}^M (b_k  - b_{k-1} ) \sum_{i= k}^M  a_i \big(\sum_{j=i}^M a_j\big)^s + \bigg(\sum_{k=N}^M  a_k \big(\sum_{i=k}^M a_i\big)^s\bigg) b_N\\
 & \approx \sum_{k=N+1}^M (b_k  - b_{k-1} ) \bigg(\sum_{i= k}^M a_i\bigg)^{s+1} + \bigg(\sum_{k=N}^M a_k\bigg)^{s+1} b_N,
 \end{align*}
in which the multiplicative constants depend only on $s$. By letting $N$ go to $-\infty$, we obtain \eqref{difference-u}.

Second, assume that $\lim_{k\to-\infty}b_k = 0$. It follows that $b_k=\sum_{i=-\infty}^k(b_i - b_{i-1})$ for every $k\in\Z$, $k\leq M$. Therefore, we have that
\begin{equation*}
\sum_{k=-\infty}^M c_k b_k =  \sum_{k=-\infty}^M (b_k  - b_{k-1}) \sum_{i= k}^M c_i
\end{equation*}
for every sequence $\{c_k\}_{k=-\infty}^M$ of  nonnegative numbers. By taking $c_k = a_k \big(\sum_{i=k}^M a_i\big)^s$ for $k\in\Z$, $k\leq M$, and using the power rules as above, we obtain \eqref{difference-u}.
 \end{proof}

The proof of the following lemma can be found in \cite[Proposition~2.1]{GolHeinStep} and \cite[Lemmas~3.2 - 3.4]{GoPi:03}.

\begin{lem}
Let $N \in \Z\cup\{-\infty\}$, $M \in \Z\cup\{+\infty\}$, $N<M$, $\beta >0$, and let $\{a_k\}_{k=N}^{M} $  and $\{\varrho_k\}_{k=N}^{M}$ be sequences of positive numbers. If $\{\varrho_k\}_{k=N}^{M}$ is nondecreasing, then
	\begin{equation}
		\sup_{N\leq k\leq M} \varrho_k  \sup_{k \leq i \leq M} a_i = \sup_{N\leq k\leq M} \varrho_k  a_k \label{inc-sup-sup}.
	\end{equation}
	If $\{\varrho_k\}_{k=N}^{M}$ is strongly increasing, then
	\begin{equation}
		\sum_{k=N}^{M} \varrho_k  \bigg(\sum_{i=k}^{M}a_i\bigg)^{\beta}  \approx \sum_{k=N}^{M} \varrho_k a_k^{\beta}, \label{inc-sum-sum}
	\end{equation}
	\begin{equation}
		\sum_{k=N}^{M} \varrho_k   \sup_{k \leq i \leq M} a_i  \approx \sum_{k=N}^{M} \varrho_k  a_k, \label{inc-sum-sup}
	\end{equation}
	and
	\begin{equation}
		\sup_{N\leq k \leq {M}} \varrho_k \bigg(\sum_{i=k}^{M}a_i\bigg)^{\beta} \approx \sup_{N\leq k \leq {M}} \varrho_k a_k^{\beta}. \label{inc-sup-sum}
	\end{equation}
	Moreover, the multiplicative constants depend only on $\inf\big\{ \frac{\varrho_{k+1}}{\varrho_k},\quad N\leq k < M\big\}$ and $\beta$.
\end{lem}
We should note that in \cite{GoPi:03}, \eqref{inc-sup-sup} is formulated when $\{\varrho_k\}_{k=N}^{M}$ is a strongly increasing sequence (and $N=-\infty)$; however, the result is a consequence of the interchanging suprema and holds true even when $\{\varrho_k\}_{k=N}^{M}$ is just nondecreasing.

\begin{defi}
Let $G$ be a positive continuous increasing function on $(a,b)$ such that $\lim_{t\to a^+}G(t)=0$. Define
\begin{equation*}
    M=\inf\{k\in\mathbb Z: G(t)\le2^k\ \text{for every $t\in(a,b)$}\}
\end{equation*}
(if the set is empty, then $M=\infty$). An increasing sequence $\{x_k\}_{k=-\infty}^M\subset (a,b]$ such that $(a,b)\subset\bigcup_{k=-\infty}^{M}[x_{k-1},x_k]$ is said to be the \textit{discretizing sequence} of $G$ if it satisfies $G(x_k) = 2^{k}$ for every $k < M$.
\end{defi}

We note that if $\lim_{t\to b^-}G(t)<\infty$, then $M < \infty$ and $x_{M}= b$, while, if $\lim_{t\to b^-}G(t)=\infty$, then $M=\infty$ and $\lim_{k\rightarrow \infty} x_k= b$. Furthermore, if $M<\infty$, then $2^{M-1} \leq G(t)\leq 2^M$ for every $t\in[x_{M-1},b)$. Note that the discretizing sequence (as defined above) is unique, and so the definite article is justified. Finally, the Darboux property of continuous functions implies that the discretizing sequence exists for every $G$ as above.

For a locally integrable nonnegative function $w$ on $[a,b)$, we will use the notation
\begin{equation*}
    W(t) = \int_a^t w(s)  \,ds, \quad t\in [a,b].
\end{equation*}
Note that the discretizing sequence for $W$ exists when $w$ is such a function.

Recall that if $M=\infty$, then $M-1$ is interpreted as $\infty$.
\begin{lem}
Let $\alpha \geq 0$. Assume that $w$ is a weight on $(a, b)$, $\{x_k\}_{k=-\infty}^{M}$ is the discretizing sequence of $W$ and $h$ is a nonnegative nonincreasing function on $(a, b)$. Then

\begin{equation}\label{int.equiv}
    \int_{a}^b W(t)^{\alpha} w(t) h(t)  \,dt \approx \sum_{k=-\infty}^{M-1} 2^{k(\alpha +1)} h(x_k)
\end{equation}
holds, in which the equivalence constants depend only on $\alpha$.
\end{lem}

\begin{proof}
The monotonicity of $h$ and properties of the discretizing sequence $\{x_k\}_{k=-\infty}^{M}$ give
\begin{align*}
\int_a^b h(t) W(t)^{\alpha} w(t) \,dt &= \sum_{k=-\infty}^{M-1} \int_{x_{k}}^{x_{k+1}} h(t) W(t)^{\alpha} w(t) \,dt \\
& \lesssim \sum_{k=-\infty}^{M-1}  h(x_k) \int_{x_k}^{x_{k+1}} d\big[W(t)^{\alpha+1} \big] \\
& \approx \sum_{k=-\infty}^{M-1}  2^{k(\alpha+1)}  h(x_k),
\end{align*}
and, conversely,
\begin{align*}
\int_a^b h(t) W(t)^{\alpha} w(t) \,dt & \geq \sum_{k=-\infty}^{M-1} \int_{x_{k-1}}^{x_k} h(t) W(t)^{\alpha} w(t) \,dt \\
& \gtrsim \sum_{k=-\infty}^{M-1}  h(x_k) \int_{x_{k-1}}^{x_k}  d\big[W(t)^{\alpha+1} \big] \\
& \approx \sum_{k=-\infty}^{M-1}  2^{k(\alpha+1)}  h(x_k).
\end{align*}
Therefore, the statement follows.
\end{proof}

Having established necessary background material, we can now take the first step towards an effective discretization of the inequality~\eqref{main-iterated-intro}.
\begin{prop}\label{P:dyadic}
Let $0 < p \leq 1$, $0 < q, r < \infty$ and let $u, v, w$ be weights on $(a,b)$.  Assume that $\{x_k\}_{k=-\infty}^{M}$ is the discretizing sequence of $W$. Then there exists a positive constant $C$ such that the inequality \eqref{main-iterated-intro} holds for all nonnegative measurable $f$ on $(a, b)$ if and only if there exist positive constants $C'$ and $C''$ such that
\begin{equation}\label{main-iterated-1}
    \bigg( \sum_{k=-\infty}^{M-1}  2^k \bigg( \int_{x_k}^{x_{k+1}} \bigg(\int_{x_k}^t f^p(s) v(s) \,ds \bigg)^\frac{q}{p} u(t) \,dt \bigg)^{\frac{r}{q}}  \bigg)^{\frac{1}{r}} \leq C'  \sum_{k=-\infty}^{M-1} \int_{x_k}^{x_{k+1}} f(t) \,dt
\end{equation}	
and
\begin{equation}\label{main-iterated-2}
    \bigg( \sum_{k=-\infty}^{M-1}  2^k   \bigg( \int_a^{x_k} f^p(t) v(t) \,dt \bigg)^\frac{r}{p} \bigg( \int_{x_k}^b  u(t) \,dt \bigg)^{\frac{r}{q}}  \bigg)^{\frac{1}{r}}  \leq C'' \sum_{k=-\infty}^{M-1} \int_{x_k}^{x_{k+1}} f(t) \,dt
\end{equation}	
for all nonnegative measurable functions $f$ on $(a, b)$. Moreover, the best constants $C$, $C'$ and $C''$ in \eqref{main-iterated-intro}, \eqref{main-iterated-1} and \eqref{main-iterated-2}, respectively, satisfy $C \approx C' + C''$.
\end{prop}	

\begin{proof}
Applying \eqref{int.equiv} with $\alpha=0$ and then using \eqref{inc-sum-sum}, we obtain
\begin{align}
    LHS\eqref{main-iterated-intro} & \approx \bigg( \sum_{k=-\infty}^{M-1}  2^k \bigg( \int_{x_k}^{x_{k+1}} \bigg(\int_a^t f^p v \bigg)^\frac{q}{p} u(t) \,dt \bigg)^{\frac{r}{q}}  \bigg)^{\frac{1}{r}} \notag \\
    & \approx \bigg( \sum_{k=-\infty}^{M-1}  2^k \bigg( \int_{x_k}^{x_{k+1}} \bigg(\int_{x_k}^t f^p v \bigg)^\frac{q}{p} u(t) \,dt \bigg)^{\frac{r}{q}}  \bigg)^{\frac{1}{r}} \notag \\
    	& \hspace{0.5 cm} + \bigg( \sum_{k=-\infty}^{M-1}  2^k   \bigg(\int_a^{x_k} f^p v \bigg)^\frac{r}{p} \bigg( \int_{x_k}^{x_{k+1}} u \bigg)^{\frac{r}{q}}  \bigg)^{\frac{1}{r}} \notag\\
     & \approx \bigg( \sum_{k=-\infty}^{M-1}  2^k \bigg( \int_{x_k}^{x_{k+1}} \bigg(\int_{x_k}^t f^p v \bigg)^\frac{q}{p} u(t) \,dt \bigg)^{\frac{r}{q}}  \bigg)^{\frac{1}{r}} \notag \\
     & \hspace{0.5 cm} + \bigg( \sum_{k=-\infty}^{M-1}  2^k   \bigg(\int_a^{x_k} f^p v \bigg)^\frac{r}{p} \bigg( \int_{x_k}^b  u \bigg)^{\frac{r}{q}}  \bigg)^{\frac{1}{r}}. \notag
\end{align}	
In the last equivalence we have used the fact that either $\int_a^{x_{M-1}} f^p v = 0$, where $x_{M-1}$ is to be interpreted as $b$ if $M=\infty$, or there is $N\in\Z\cup\{-\infty\}$, $N\leq M-1$, such that $\int_a^{x_k} f^p v=0$ for every $k<N$ and $\bigg\{2^k \bigg(\int_a^{x_k} f^p v \bigg)^\frac{r}{p}\bigg\}_{k=N}^{M-1}$ is a strongly increasing sequence (unless $N=M-1<\infty$, which is a trivial case). The assertion follows.
\end{proof}

The next step is based on saturation of Hardy inequalities and embeddings of weighted Lebesgue spaces on the intervals determined by a discretizing sequence.

\begin{prop}\label{P:2-inequalities}
Let $0 < p \leq 1$, $0 < q, r < \infty$ and let $u, v, w$ be weights on $(a,b)$.  Assume that $\{x_k\}_{k=-\infty}^{M}$ is the discretizing sequence of $W$. For every $k\in\Z$, $k\leq M$, we denote
\begin{equation}\label{A(a,b)}
    A_k := \sup_{g\in  \Mpl(x_{k-1},x_k)} \frac{\bigg(\int_{x_{k-1}}^{x_k} g(t)^p v(t) \,dt \bigg)^{\frac{1}{p}}} {\int_{x_{k-1}}^{x_k} g(t) \,dt}
\end{equation}
and
\begin{equation}\label{B(a,b)}
    B_k := \sup_{h\in  \Mpl(x_{k-1},x_{k})} \frac{\bigg(\int_{x_{k-1}}^{x_{k}} \bigg(\int_{x_{k-1}}^t h(s)^p v(s) \,ds \bigg)^{\frac{q}{p}} u(t) \,dt \bigg)^{\frac{1}{q}}}{\int_{x_{k-1}}^{x_{k}}  h(t) \,dt}.
\end{equation}
Then there exists a positive constant $C$ such that the inequality \eqref{main-iterated-intro} holds for all nonnegative measurable functions $f$ on $(a, b)$ if and only if there exist positive constants $C', C''$ such that the inequalities
\begin{equation}\label{Bk-inequality}
\bigg(\sum_{k=-\infty}^{M-1} 2^k a_k^r B_{k+1}^r \bigg)^{\frac{1}{r}} \leq C'  \sum_{k=-\infty}^{M-1} a_k
\end{equation}
and
\begin{equation}\label{Ak-inequality}
\bigg(\sum_{k=-\infty}^{M-1} 2^k \bigg(\int_{x_k}^b u(t) \,dt \bigg)^{\frac{r}{q}} \bigg(\sum_{j=-\infty}^{k} a_j^p A_j^p \bigg)^\frac{r}{p}\bigg)^{\frac{1}{r}} \leq C'' \sum_{k=-\infty}^{M-1} a_k
\end{equation}	
hold for every sequence $\{a_k\}_{k=-\infty}^{M-1}$ of nonnegative numbers. Moreover, the best constants $C$, $C'$ and $C''$ in \eqref{main-iterated-intro}, \eqref{Bk-inequality} and \eqref{Ak-inequality}, respectively, satisfy $C \approx C' + C''$.
\end{prop}	

\begin{proof}
Assume that \eqref{main-iterated-1} holds. By~\eqref{B(a,b)}, there exist nonnegative measurable functions $h_k$, $k \leq M-1$, on $(a, b)$ such that $\supp h_k \subset [x_k, x_{k+1}]$, $\int_{x_{k}}^{x_{k+1}} h_k = 1$, and
\begin{equation*}
    \bigg(\int_{x_k}^{x_{k+1}} \bigg(\int_{x_k}^t h_k^p v \bigg)^{\frac{q}{p}} u(t) \,dt\bigg)^{\frac{1}{q}} \gtrsim B_{k+1}.
\end{equation*}
Thus, given $\{a_m\}_{m=-\infty}^{M-1}$ and inserting $f = \sum_{m=-\infty}^{M-1}  a_m h_m$   into \eqref{main-iterated-1}, we get \eqref{Bk-inequality}. Conversely, \eqref{main-iterated-1} follows at once from~\eqref{Bk-inequality} on setting  $a_k = \int_{x_k}^{x_{k+1}} f$ for $k\in(-\infty,M-1)$.

Similarly, by~\eqref{A(a,b)}, there exist nonnegative measurable functions $g_k$, $k \leq M-1$, on $(a, b)$ such that $\supp g_k \subset [x_{k-1}, x_k]$, $\int_{x_{k-1}}^{x_{k}} g_k = 1$, and
\begin{equation*}
    \bigg(\int_{x_{k-1}}^{x_{k}} g_k^p v \bigg)^{\frac{1}{p}} \gtrsim A_k.
\end{equation*}
Thus, given $\{a_m\}_{m=-\infty}^{M-1}$ and inserting $f = \sum_{m=-\infty}^{M-1}  a_m g_m$ into \eqref{main-iterated-2}, \eqref{Ak-inequality} follows. Conversely, inserting  $a_k = \int_{x_{k-1}}^{x_k} f$ in \eqref{Ak-inequality} gives \eqref{main-iterated-2}.

The assertion now directly follows from Proposition~\ref{P:dyadic}.
\end{proof}

\section{Proofs}\label{S:proofs}

We begin this section with a theorem of auxiliary nature, albeit interesting on its own, which yields a discrete characterization of the inequality in question. We will then use it as the last step towards the proof of Theorem~A.

\begin{thm}\label{T:disc-p<1}
Let $0 < p \leq 1$, $0 < q, r < \infty$ and let $u, v, w$ be weights on $(a,b)$. Let $\{x_k\}_{k=-\infty}^{M}$ be the discretizing sequence of $W$. Then there exists a constant $C>0$ such that the inequality \eqref{main-iterated-intro} holds for all nonnegative measurable functions $f$ on $(a, b)$ if and only if one of the following conditions is satisfied:
	
\rm(i) $ 1 \le r $, $1\leq q$,
\begin{align*}
A_1^* &:=\sup_{ k \leq M-1} 2^{\frac{k}{r}}  \esup_{t \in (x_{k}, x_{k+1})} \bigg(\int_t^{x_{k+1}} u(s) \,ds \bigg)^{\frac{1}{q}} V_p(x_k, t)  < \infty \notag\\
\intertext{and}
B_1^* &:=  \sup_{k \leq M-1} \bigg(\sum_{i=k}^{M-1} 2^{i} \bigg(\int_{x_i}^{b} u(t) \,dt \bigg)^{\frac{r}{q}} \bigg)^{\frac{1}{r}} V_p(a, x_k) < \infty;
\end{align*}

\rm(ii) $  1 \le r $, $q<1$, $B_1^* < \infty$ and
\begin{equation*}
A_2^* :=\sup_{k \leq {M-1}} 2^{\frac{k}{r}} \bigg(\int_{x_k}^{x_{k+1}} \bigg(\int_t^{x_{k+1}} u(s) \,ds \bigg)^{\frac{q}{1-q}} u(t) V_p(x_k, t)^{\frac{q}{1-q}}  \,dt \bigg)^{\frac{1-q}{q}} < \infty;
\end{equation*}

\rm(iii) $r < 1$, $1 \leq q$,
\begin{align*}
A_3^* &:=\bigg(\sum_{k=-\infty}^{M-1} 2^{\frac{k}{1-r}} \esup_{t \in (x_k, x_{k+1})} \bigg(\int_t^{x_{k+1}} u(s) \,ds\bigg)^{\frac{r}{q(1-r)}} V_p(x_k, t)^{\frac{r}{1-r}} \bigg)^{\frac{1-r}{r}}< \infty, \notag \\
\intertext{and}
B_2^* &:= \bigg( \sum_{k=-\infty}^{{M-1}} 2^{k} \bigg(\int_{x_k}^b u(t) \,dt \bigg)^{\frac{r}{q}} \bigg(\sum_{i=k }^{{M-1}} 2^i \bigg(\int_{x_i}^{b} u(t) \,dt \bigg)^{\frac{r}{q}}  \bigg)^{\frac{r}{1-r}} 	V_p(a, x_k)^{\frac{r}{1-r}} \bigg)^{\frac{1-r}{r}} < \infty;
\end{align*}

\rm(iv) $r < 1$, $q < 1$, $B_2^* < \infty$ and
\begin{equation*}
A_4^* := \bigg( \sum_{k=-\infty}^{{M-1}} 2^{\frac{k}{1-r}} \bigg( \int_{x_k}^{x_{k+1}} \bigg(\int_{t}^{x_{k+1}} u(s) \,ds \bigg)^{\frac{q}{1-q}} u(t) V_p(x_k, t)^{\frac{q}{1-q}} \,dt \bigg)^{\frac{r(1-q)}{q(1-r)}} \bigg)^{\frac{1-r}{r}} < \infty.
\end{equation*}

Moreover, the best constant $C$ in the inequality \eqref{main-iterated-intro} satisfies
\begin{equation*}
    C\approx
    \begin{cases}
        A_1^* + B_1^* &\text{in the case (i),}
            \\
        A_2^* + B_1^* &\text{in the case (ii),}
            \\
        A_3^* + B_2^* &\text{in the case (iii),}
            \\
        A_4^* + B_2^* &\text{in the case (iv).}
    \end{cases}
\end{equation*}
\end{thm}

\begin{proof}
It follows from Proposition~\ref{P:2-inequalities} that the best constant $C$ in \eqref{main-iterated-intro} satisfies $C \approx C' + C''$, where $C'$ and $C''$ are the best constants in \eqref{Bk-inequality} and \eqref{Ak-inequality}. Next, we obtain an appropriate characterization of $C'$ by combining a discrete version of the Landau resonance theorem (cf.~e.g.~\cite[Proposition~4.1]{GoPi:03}) with the classical Hardy inequality. Finally, an appropriate two-sided estimate of $C''$ can be obtained by combining the known characterization of a discrete Hardy inequality (cf.~e.g.~\cite[Theorem~1]{Be:91} or~\cite[Theorem~9.2]{Gr:98})) with the classical duality expression of the norm in a weighted Lebesgue space.
\end{proof}

\begin{proof}[Proof of Theorem A]
First of all, note that the optimal constant $C$ in \eqref{main-iterated-intro} is equal to $\infty$ if there is $t_0\in(a,b)$ such that $W(t_0)=\infty$, and so is RHS\eqref{thm:optimal_constant}; hence the theorem is trivially true in this pathological case. Therefore, we may assume that $W(t)<\infty$ for every $t\in(a,b)$. Let $\{x_k\}_{k=-\infty}^M$, where $M\in\Z\cup\{\infty\}$, be the discretizing sequence of $W$.

\rm(i) Let $ p\leq  1 \le r $, $1\leq q$. We have from Theorem~\ref{T:disc-p<1}(i) that $C \approx A_1^* + B_1^*$.  Define
\begin{equation*}
\tilde{A}_1 :=\sup_{k \leq M-1} 2^{\frac{k}{r}}  \esup_{t \in (x_{k}, b)} \bigg(\int_t^b u \bigg)^{\frac{1}{q}} V_p(a, t).
\end{equation*}
We will first show that $A_1^* + B_1^* \approx \tilde{A}_1 + B_1^*$. Since obviously $A_1^* \leq \tilde{A}_1$, it is enough to prove that $\tilde{A}_1 \lesssim A_1^* + B_1^*$. Using \eqref{inc-sup-sup}, we obtain
\begin{align*}
\tilde{A}_1 & = \sup_{k \leq M-1} 2^{\frac{k}{r}} \sup_{k \leq i \leq M-1} \esup_{t \in (x_{i}, x_{i+1})} \bigg(\int_t^b u \bigg)^{\frac{1}{q}} V_p(a, t) =\sup_{k \leq M-1} 2^{\frac{k}{r}}  \esup_{t \in (x_{k}, x_{k+1})} \bigg(\int_t^b u \bigg)^{\frac{1}{q}} V_p(a, t) \\
& \approx \sup_{k \leq M-1} 2^{\frac{k}{r}}  \esup_{t \in (x_{k}, x_{k+1})} \bigg(\int_t^{x_{k+1}} u \bigg)^{\frac{1}{q}} V_p(a, t)  + \sup_{k \leq M-2} 2^{\frac{k}{r}} \bigg(\int_{x_{k+1}}^b u \bigg)^{\frac{1}{q}} V_p(a, x_{k+1}).
\end{align*}
Since
\begin{equation}\label{V-cut}
    V_p(a,t)\approx V_p(a,x_k)+V_p(x_k,t) \quad\text{for every $t\in(x_k,x_{k+1})$},
\end{equation}
we in fact have
\begin{align*}
\tilde{A}_1 &  \approx A_1^* + \sup_{k \leq M-1} 2^{\frac{k}{r}}  \bigg(\int_{x_k}^{x_{k+1}} u \bigg)^{\frac{1}{q}} V_p(a, x_k)
+  \sup_{k \leq M-2} 2^{\frac{k}{r}} \bigg(\int_{x_{k+1}}^b u \bigg)^{\frac{1}{q}} V_p(a, x_{k+1})
    \\
&\lesssim A_1^* + \sup_{k \leq M-1} 2^{\frac{k}{r}}  \bigg(\int_{x_k}^{b} u \bigg)^{\frac{1}{q}} V_p(a, x_k)  \le A_1^* + B_1^*,
\end{align*}
establishing the claim.

Next, we will show that $C_1+ C_2 \approx \tilde{A}_1 + B_1^* $. Observe first that,
\begin{align*}
C_1 &= \sup_{k \leq M-1} \sup_{ t \in (x_k, x_{k+1})} \bigg(\int_a^t w \bigg)^{\frac{1}{r}} \esup_{s \in (t,b)} \bigg(\int_s^b u \bigg)^{\frac{1}{q}} V_p(a, s) \approx \tilde{A}_1.
\end{align*}

On the other hand, fixing $k \in \Z$, $k <M$, we have that
\begin{align}\label{int-uw-sum-upper}
\sum_{i=k}^{M-1} 2^i \bigg(\int_{x_i}^b u \bigg)^{\frac{r}{q}} &= 2^k \bigg(\int_{x_k}^b u \bigg)^{\frac{r}{q}} +  \sum_{i=k+1}^{M-1} 2^i \bigg(\int_{x_i}^b u \bigg)^{\frac{r}{q}} \notag\\
& \approx 2^k \bigg(\int_{x_k}^b u \bigg)^{\frac{r}{q}} +  \sum_{i=k+1}^{M-1} \bigg( \int_{x_{i-1}}^{x_i} w \bigg) \bigg(\int_{x_i}^b u \bigg)^{\frac{r}{q}} \notag \\
& \leq 2^k \bigg(\int_{x_k}^b u \bigg)^{\frac{r}{q}} +  \int_{x_{k}}^b w(t)  \bigg(\int_t^b u \bigg)^{\frac{r}{q}} \,dt
\end{align}
with equivalence constants independent of $k$. Hence, in view of \eqref{int-uw-sum-upper}, we have
\begin{align}
B_1^*& \lesssim  \sup_{k \leq M-1} \bigg(\int_{x_k}^b w(t) \bigg(\int_t^b u \bigg)^{\frac{r}{q}} \,dt \bigg)^{\frac{1}{r}}V_p(a, x_k) + \sup_{k \leq M-1} 2^\frac{k}{r}  \bigg(\int_{x_k}^b  u \bigg)^{\frac{r}{q}} V_p(a, x_k) \notag\\
& \leq \sup_{k \leq M-1} \esup_{t \in (x_k, x_{k+1})} \bigg(\int_{t}^b w(s) \bigg(\int_s^b u \bigg)^{\frac{r}{q}} \,ds \bigg)^{\frac{1}{r}} V_p(a, t) +\tilde{A}_1 \notag \\
& \approx  C_2 + C_1. \label{B1*<C2}
\end{align}
Thus, we have $\tilde{A}_1 + B_1^* \lesssim C_1 + C_2$.

Conversely,
\begin{equation} \label{int-uw-sum-lower}
\int_{x_k}^b w(t) \bigg(\int_t^b u \bigg)^{\frac{r}{q}} \,dt = \sum_{i=k}^{M-1} \int_{x_i
}^{x_{i+1}} w(t) \bigg(\int_t^b u \bigg)^{\frac{r}{q}} \,dt \lesssim \sum_{i=k}^{M-1} 2^i \bigg(\int_{x_i}^b u \bigg)^{\frac{r}{q}}.
\end{equation}
Consequently,
\begin{align*}
C_2 & \approx \sup_{k \leq M-1 } \esup_{ t \in (x_k, x_{k+1})} \bigg(\int_t^{x_{k+1}} w(s) \bigg(\int_s^b u \bigg)^{\frac{r}{q}} \,ds \bigg)^{\frac{1}{r}} V_p(a, t)\\
    & \quad + \sup_{k \leq M-2 } \bigg(\int_{x_{k+1}}^b w(t) \bigg(\int_t^b u \bigg)^{\frac{r}{q}} \,dt \bigg)^{\frac{1}{r}} V_p(a, x_{k+1}).
\end{align*}
Hence, in view of \eqref{int-uw-sum-lower},~\eqref{inc-sup-sup} and the fact that $\{x_k\}_{k=-\infty}^M$ is the discretizing sequence for $W$, we obtain
\begin{align}
C_2 & \lesssim  \sup_{k \leq M-1 } 2^{\frac{k}{r}} \esup_{ t \in (x_k, x_{k+1})}  \bigg(\int_t^b u \bigg)^{\frac{1}{q}} V_p(a, t) \notag\\
    & \quad + \sup_{k \leq M-2 } \bigg(\sum_{i=k+1}^{M-1} 2^i \bigg(\int_{x_i}^{b} u \bigg)^{\frac{r}{q}} \bigg)^{\frac{1}{r}} V_p(a, x_{k+1}) \notag\\
& \lesssim \tilde{A}_1 + B_1^*.\label{C2<A1+B1*}
\end{align}
Consequently, we arrive at $C \approx C_1 + C_2$.

\rm(ii) Let $ p\leq 1 \le r $, $q<1$. Using Theorem~\ref{T:disc-p<1}(ii), we have that $C \approx A_2^* + B_1^*$. Let us show that $A_2^* + B_1^* \approx \tilde{A}_2 + B_1^*$, where
\begin{equation*}
\tilde{A}_2 := \sup_{k \leq {M-1}} 2^{\frac{k}{r}} \bigg(\int_{x_k}^{b} \bigg(\int_t^{b} u \bigg)^{\frac{q}{1-q}} u(t) V_p(a, t)^{\frac{q}{1-q}}  \,dt \bigg)^{\frac{1-q}{q}}.
\end{equation*}
It is easy to see that $A_2^* \lesssim \tilde{A}_2$. On the other hand, using \eqref{inc-sup-sum}, we have
\begin{align*}
\tilde{A}_2&= \sup_{k \leq {M-1}} 2^{\frac{k}{r}} \bigg(\sum_{i=k}^{M-1}\int_{x_i}^{x_{i+1}} \bigg(\int_t^{b} u \bigg)^{\frac{q}{1-q}}u(t) V_p(a, t)^{\frac{q}{1-q}} \,dt  \bigg)^{\frac{1-q}{q}}\\	& \approx \sup_{k \leq {M-1}} 2^{\frac{k}{r}} \bigg(\int_{x_k}^{x_{k+1}} \bigg(\int_t^{b} u \bigg)^{\frac{q}{1-q}} u(t) V_p(a, t)^{\frac{q}{1-q}}  \,dt \bigg)^{\frac{1-q}{q}}.
\end{align*}
Furthermore, for each $k \leq M-1$, we have that
\begin{align}\label{IBP}
& \bigg(\int_{x_k}^{x_{k+1}} \bigg(\int_t^{b} u \bigg)^{\frac{q}{1-q}} u(t) V_p(a, t)^{\frac{q}{1-q}}  \,dt \bigg)^{\frac{1-q}{q}} \notag\\
	&\hspace{2cm} \lesssim \bigg(\int_{x_k}^{x_{k+1}} \bigg(\int_t^{b} u \bigg)^{\frac{1}{1-q}} d\big[V_p(a, t)^{\frac{q}{1-q}}\big] \bigg)^{\frac{1-q}{q}} + \lim_{t \rightarrow x_k+} \bigg(\int_t^{b} u \bigg)^{\frac{1}{q}} V_p(a, t).
\end{align}
Indeed, by integrating by parts, it is clear that \eqref{IBP} holds for each $k \in \Z$, $k < M-1$, whether $M=\infty$ or $M<\infty$. The remaining case when $M<\infty$ and $k = M-1$ requires more explanation. We may assume that $\max\{A_2^*, B_1^*\}< \infty$; consequently
\begin{align*}
& \bigg(\int_{x_{M-1}}^{b} \bigg(\int_t^{b} u \bigg)^{\frac{q}{1-q}} u(t) V_p(a, t)^{\frac{q}{1-q}}  \,dt \bigg)^{\frac{1-q}{q}} < \infty.
\end{align*}
Thus, for each $x \in (x_{M-1}, b)$,
\begin{align*}
V_p(a,x) \bigg(\int_x^{b} u \bigg)^{\frac{1}{q}} \lesssim \bigg(\int_{x}^{b} \bigg(\int_t^{b} u \bigg)^{\frac{q}{1-q}} u(t) V_p(a, t)^{\frac{q}{1-q}}  \,dt \bigg)^{\frac{1-q}{q}}
\end{align*}
holds, whence we conclude that
\begin{align*}
	\lim_{x \rightarrow b-} V_p(a,x) \bigg(\int_x^{b} u \bigg)^{\frac{1}{q}} = 0.
\end{align*}
Hence, \eqref{IBP} holds.
Additionally, observe that
\begin{align} \label{limit<supremum}
	\lim_{t \rightarrow x_k+} \bigg(\int_t^{b} u \bigg)^{\frac{1}{q}} V_p(a, t) \leq \esup_{t \in (x_k, x_{k+1})} \bigg(\int_t^{b} u \bigg)^{\frac{1}{q}} V_p(a, t).
\end{align}
Then, in view of \eqref{IBP} and \eqref{limit<supremum},

\begin{align*}
\tilde{A}_2&\lesssim \sup_{k \leq {M-1}} 2^{\frac{k}{r}} \bigg(\int_{x_k}^{x_{k+1}} \bigg(\int_t^{b} u \bigg)^{\frac{1}{1-q}} d\big[V_p(a, t)^{\frac{q}{1-q}}\big] \bigg)^{\frac{1-q}{q}}\\
&\quad + \sup_{k \leq {M-1}} 2^{\frac{k}{r}} \esup_{t \in (x_k, x_{k+1})} \bigg(\int_t^{b} u \bigg)^{\frac{1}{q}} V_p(a, t) \\
& \lesssim \sup_{k \leq {M-1}} 2^{\frac{k}{r}} \bigg(\int_{x_k}^{x_{k+1}} \bigg(\int_t^{x_{k+1}} u \bigg)^{\frac{1}{1-q}} d\big[V_p(a, t)^{\frac{q}{1-q}}\big] \bigg)^{\frac{1-q}{q}}\\
&\quad + \sup_{k \leq {M-1}} 2^{\frac{k}{r}} \esup_{t \in (x_k, x_{k+1})} \bigg(\int_t^{x_{k+1}} u \bigg)^{\frac{1}{q}} V_p(a, t) \\
& \quad + \sup_{k \leq {M-2}} 2^{\frac{k}{r}} \bigg(\int_{x_{k+1}}^b u \bigg)^{\frac{1}{q}} V_p(a, x_{k+1})\\
& =: \tilde{A}_{2,1} + \tilde{A}_{2,2} + \tilde{A}_{2,3}.
\end{align*}
We shall now establish appropriate upper estimates for $\tilde{A}_{2,1}$, $\tilde{A}_{2,2}$ and $\tilde{A}_{2,3}$. Note that, \eqref{V-cut} yields that
\begin{equation}\label{lower-A2*_+_B1*}
\begin{aligned}
&\sup_{k \leq {M-1}} 2^{\frac{k}{r}} \bigg(\int_{x_k}^{x_{k+1}} \bigg(\int_t^{x_{k+1}} u \bigg)^{\frac{q}{1-q}} u(t) \, V_p(a, t)^{\frac{q}{1-q}}  \,dt \bigg)^{\frac{1-q}{q}} \\
&\hspace{0.5 cm}\approx A_2^* +  \sup_{k \leq {M-1}} 2^{\frac{k}{r}} \bigg(\int_{x_k}^{x_{k+1}} u \bigg)^{\frac{1}{q}} V_p(a, x_{k})\\
&\hspace{0.5 cm} \leq A_2^* + B_1^*.
\end{aligned}
\end{equation}
Since integration by parts gives
\begin{align*}
\tilde{A}_{2,1} \lesssim  \sup_{k \leq {M-1}} 2^{\frac{k}{r}} \bigg(\int_{x_k}^{x_{k+1}} \bigg(\int_t^{x_{k+1}} u \bigg)^{\frac{q}{1-q}} u(t) \, V_p(a, t)^{\frac{q}{1-q}}  \,dt \bigg)^{\frac{1-q}{q}},
\end{align*}
it follows that $\tilde{A}_{2,1}\lesssim A_2^* + B_1^*$. Furthermore, note that
\begin{align}
\esup_{t \in (x, y)} \bigg(\int_t^y u \bigg)^{\frac{1}{q}} V_p(a, t) & \approx \esup_{t \in (x, y)} \bigg(\int_t^y \bigg(\int_s^y u \bigg)^{\frac{q}{1-q}} u(s) \,ds \bigg)^{\frac{1-q}{q}} V_p(a, t) \notag\\
&  \leq  \bigg(\int_x^y \bigg(\int_t^y u \bigg)^{\frac{q}{1-q}} u(t) V_p(a, t)^{\frac{q}{1-q}} \,dt \bigg)^{\frac{1-q}{q}}. \label{upper-A1-A2}
\end{align}
Thus, applying \eqref{upper-A1-A2} and \eqref{lower-A2*_+_B1*}, we obtain that
\begin{equation*}
\tilde{A}_{2,2} \lesssim \sup_{k\leq M-1} 2^{\frac{k}{r}} \bigg(\int_{x_k}^{x_{k+1}} \bigg(\int_t^{x_{k+1}} u \bigg)^{\frac{q}{1-q}} u(t) V_p(a, t)^{\frac{q}{1-q}} \,dt \bigg)^{\frac{1-q}{q}} \lesssim A_2^* + B_1^*.
\end{equation*}
Finally, it is clear that $\tilde{A}_{2,3} \lesssim B_1^*$. Combining these estimates we arrive at $C \approx \tilde{A}_2 + B_1^*$.

Next, we will prove that $\tilde{A}_2 + B_1^* \approx C_2 + C_3$. We have from \eqref{B1*<C2} that $B_1^* \lesssim C_2$. Moreover,
\begin{align*}
C_3 = \sup_{k \leq M-1} \sup_{t\in (x_k, x_{k+1})} \bigg(\int_a^t w \bigg)^{\frac{1}{r}}  \bigg(\int_{t}^{b} \bigg(\int_s^b u \bigg)^{\frac{q}{1-q}} u(s) V_p(a, s)^{\frac{q}{1-q}} \,ds \bigg)^{\frac{1-q}{q}} \approx \tilde{A}_2.
\end{align*}
Therefore, the proof will be complete once we show that $C_2 \lesssim \tilde{A}_2 + B_1^*$.

Applying \eqref{upper-A1-A2}, we plainly have that $\tilde{A}_1 \lesssim \tilde{A}_2$. Finally, using \eqref{C2<A1+B1*}, we obtain  $C_2 \lesssim
\tilde{A}_1 + B_1^* \lesssim \tilde{A}_2 + B_1^*$. Hence the proof is complete in this case.

\rm(iii) Let $p \leq 1 $, $r < 1$, $1 \leq q$, then we have from Theorem~\ref{T:disc-p<1}(iii) that $C \approx A_3^* + B_2^*$. First, we will show that $C \approx \tilde{A}_3 + B_2^*$, where
\begin{equation*}
\tilde{A}_3 :=\bigg(\sum_{k=-\infty}^{M-1} 2^{\frac{k}{1-r}} \esup_{t \in (x_k, b)} \bigg(\int_t^{b} u\bigg)^{\frac{r}{q(1-r)}} V_p(a, t)^{\frac{r}{1-r}} \bigg)^{\frac{1-r}{r}}.
\end{equation*}
It is clear that $A_3^* \leq \tilde{A}_3$. Moreover, \eqref{inc-sum-sup} together with \eqref{V-cut} yields
\begin{align*}
\tilde{A}_3 & = \bigg(\sum_{k=-\infty}^{M-1} 2^{\frac{k}{1-r}} \sup_{k \leq i \leq M-1} \esup_{t \in (x_i, x_{i+1})} \bigg(\int_t^{b} u\bigg)^{\frac{r}{q(1-r)}} V_p(a, t)^{\frac{r}{1-r}} \bigg)^{\frac{1-r}{r}} \\
& \approx \bigg(\sum_{k=-\infty}^{M-1} 2^{\frac{k}{1-r}} \esup_{t \in (x_k, x_{k+1})} \bigg(\int_t^{b} u\bigg)^{\frac{r}{q(1-r)}} V_p(a, t)^{\frac{r}{1-r}} \bigg)^{\frac{1-r}{r}} \\
& \approx \bigg(\sum_{k=-\infty}^{M-1} 2^{\frac{k}{1-r}} \esup_{t \in (x_k, x_{k+1})} \bigg(\int_t^{x_{k+1}} u\bigg)^{\frac{r}{q(1-r)}} V_p(a, t)^{\frac{r}{1-r}} \bigg)^{\frac{1-r}{r}} \\
    & \quad + \bigg(\sum_{k=-\infty}^{M-2} 2^{\frac{k}{1-r}} \bigg(\int_{x_{k+1}}^{b} u \bigg)^{\frac{r}{q(1-r)}}  V_p(a, x_{k+1})^{\frac{r}{1-r}} \bigg)^{\frac{1-r}{r}} \\
& \lesssim A_3^* + \bigg(\sum_{k=-\infty}^{M-1} 2^{\frac{k}{1-r}} \bigg(\int_{x_k}^b u\bigg)^{\frac{r}{q(1-r)}} V_p(a, x_k)^{\frac{r}{1-r}} \bigg)^{\frac{1-r}{r}} \\
& \leq A_3^* + B_2^*.
\end{align*}

Next, we will show that $\tilde{A}_3 + B_2^* \approx C_4 + C_5$. We will find equivalent formulations for $B_2^*$. Using \eqref{difference-u} for $a_k = 2^k \bigg(\int_{x_k}^{b} u \bigg)^{\frac{r}{q}}$, $b_k = V_p(a, x_k)^{\frac{r}{1-r}}$ and $s=\frac{r}{1-r}$, we get
\begin{align*}
B_2^*  &\approx  \bigg( \sum_{k=-\infty}^{{M-1}} \bigg(\sum_{i=k}^{{M-1}} 2^i \bigg(\int_{x_i}^{b} u \bigg)^{\frac{r}{q}}  \bigg)^{\frac{1}{1-r}} \big[V_p(a, x_k)^{\frac{r}{1-r}} - V_p(a, x_{k-1})^{\frac{r}{1-r}} \big]
 \bigg)^{\frac{1-r}{r}} \notag\\
 &  \hspace{2cm} +  \bigg(\sum_{i=-\infty}^{{M-1}} 2^i \bigg(\int_{x_i}^{b} u \bigg)^{\frac{r}{q}}\bigg)^{\frac{1}{r}}  \lim_{k\rightarrow -\infty} V_p(a, x_k).
\end{align*}
Applying  \eqref{int-uw-sum-upper} and \eqref{int.equiv} with $\alpha = 0$, we have that
\begin{align*}
B_2^* &\lesssim  \bigg( \sum_{k=-\infty}^{{M-1}} \bigg(\int_{x_k}^b w(t) \bigg(\int_{t}^b u \bigg)^{\frac{r}{q}} \,dt \bigg)^{\frac{1}{1-r}}  \big[V_p(a, x_k)^{\frac{r}{1-r}} - V_p(a, x_{k-1})^{\frac{r}{1-r}} \big]  \bigg)^{\frac{1-r}{r}} \\
    & \hspace{0.5cm} + \bigg( \sum_{k=-\infty}^{{M-1}} 2^{\frac{k}{1-r}}  \bigg(\int_{x_k}^b u \bigg)^{\frac{r}{q(1-r)}} V_p(a, x_k)^{\frac{r}{1-r}} \bigg)^{\frac{1-r}{r}} \\
    & \hspace{0.5cm} + \bigg(\int_a^b w(t) \bigg(\int_t^b u\bigg)^{\frac{r}{q}} \,dt \bigg)^{\frac{1}{r}} \lim_{k\rightarrow -\infty} V_p(a, x_k)\\
& \lesssim \bigg( \sum_{k=-\infty}^{{M-1}} \int_{x_k}^{b} \bigg(\int_t^b w(s) \bigg(\int_s^{b} u \bigg)^{\frac{r}{q}} \,ds  \bigg)^{\frac{r}{1-r}} w(t) \bigg(\int_t^{b} u \bigg)^{\frac{r}{q}} \,dt  \\
& \hspace{4cm} \times
\big[V_p(a, x_k)^{\frac{r}{1-r}} - V_p(a, x_{k-1})^{\frac{r}{1-r}} \big]  \bigg)^{\frac{1-r}{r}} \\
 & \hspace{0.5cm} + \bigg( \sum_{k=-\infty}^{{M-1}} 2^{\frac{k}{1-r}} \esup_{t \in (x_k, x_{k+1})} \bigg(\int_{t}^b u \bigg)^{\frac{r}{q(1-r)}} V_p(a, t)^{\frac{r}{1-r}} \bigg)^{\frac{1-r}{r}} \\
    & \hspace{0.5cm} +\bigg( \int_a^b \bigg(\int_t^b w(s) \bigg(\int_s^{b} u \bigg)^{\frac{r}{q}} \,ds  \bigg)^{\frac{r}{1-r}} w(t) \bigg(\int_t^{b} u \bigg)^{\frac{r}{q}} \,dt\bigg)^{\frac{1-r}{r}} \lim_{k\rightarrow -\infty} V_p(a, x_k). \notag
\end{align*}
Now, plugging $b_k = V_p(a, x_k)^{\frac{r}{1-r}}$ and
\begin{equation*}
    c_k = \int_{x_k}^{x_{k+1}} \bigg(\int_t^b w(s) \bigg(\int_s^{b} u \bigg)^{\frac{r}{q}} \,ds  \bigg)^{\frac{r}{1-r}} w(t) \bigg(\int_t^{b} u \bigg)^{\frac{r}{q}} \,dt
\end{equation*}
into \eqref{L:equiv.cond.:Abel}, we obtain
\begin{align}
B_2^* & \lesssim \bigg( \sum_{k=-\infty}^{{M-1}} \int_{x_k}^{x_{k+1}} \bigg(\int_t^b w(s) \bigg(\int_s^{b} u \bigg)^{\frac{r}{q}} \,ds  \bigg)^{\frac{r}{1-r}} w(t) \bigg(\int_t^{b} u \bigg)^{\frac{r}{q}}
 \,dt\, V_p(a, x_k)^{\frac{r}{1-r}} \bigg)^{\frac{1-r}{r}} \notag \\
 & \hspace{0.5cm} +  \bigg( \sum_{k=-\infty}^{{M-1}} 2^{\frac{k}{1-r}} \esup_{t \in (x_k, x_{k+1})} \bigg(\int_{t}^b u \bigg)^{\frac{r}{q(1-r)}} V_p(a, t)^{\frac{r}{1-r}} \bigg)^{\frac{1-r}{r}} \notag\\
& \leq \bigg( \sum_{k=-\infty}^{{M-1}} \int_{x_k}^{x_{k+1}} \bigg(\int_t^b w(s) \bigg(\int_s^{b} u \bigg)^{\frac{r}{q}} \,ds  \bigg)^{\frac{r}{1-r}} w(t) \bigg(\int_t^{b} u \bigg)^{\frac{r}{q}}
V_p(a, t)^{\frac{r}{1-r}} \,dt \bigg)^{\frac{1-r}{r}} + \tilde{A}_3 \notag \\
& \approx C_5 + C_4. \label{B2*<C5}
\end{align}
In the last equivalence we use $\tilde{A}_3 \approx C_4$, which easily follows from \eqref{int.equiv} with $\alpha=\frac{r}{1-r}$. Hence, we have $\tilde{A}_3 + B_2^* \lesssim C_4 + C_5$.

Conversely, integration by parts gives
\begin{align*}
C_5 & = \bigg( \sum_{k=-\infty}^{{M-1}} \int_{x_k}^{x_{k+1}} \bigg(\int_t^b w(s) \bigg(\int_s^{b} u \bigg)^{\frac{r}{q}} \,ds  \bigg)^{\frac{r}{1-r}} w(t) \bigg(\int_t^{b} u \bigg)^{\frac{r}{q}}
V_p(a, t)^{\frac{r}{1-r}} \,dt \bigg)^{\frac{1-r}{r}} \\
& \lesssim \bigg( \sum_{k=-\infty}^{{M-1}} V_p(a, x_k)^{\frac{r}{1-r}} \bigg[\bigg(\int_{x_k}^b w(t) \bigg(\int_t^{b} u \bigg)^{\frac{r}{q}} \,dt  \bigg)^{\frac{1}{1-r}}  \\
& \hspace{5cm} - \bigg(\int_{x_{k+1}}^b w(t) \bigg(\int_t^{b} u \bigg)^{\frac{r}{q}} \,dt  \bigg)^{\frac{1}{1-r}}\bigg]  \bigg)^{\frac{1-r}{r}} \\
    & \quad  + \bigg( \sum_{k=-\infty}^{{M-1}} \int_{x_k}^{x_{k+1}} \bigg(\int_t^b w(s) \bigg(\int_s^{b} u \bigg)^{\frac{r}{q}} \,ds  \bigg)^{\frac{1}{1-r}} d\big[V_p(a,t)^{\frac{r}{1-r}}\big] \bigg)^{\frac{1-r}{r}}\\
& \approx B_2^* + \bigg( \sum_{k=-\infty}^{{M-1}} \int_{x_k}^{x_{k+1}} \bigg(\int_t^{x_{k+1}} w(s) \bigg(\int_s^{b} u \bigg)^{\frac{r}{q}} \,ds  \bigg)^{\frac{1}{1-r}} d\big[V_p(a,t)^{\frac{r}{1-r}}\big] \bigg)^{\frac{1-r}{r}} \\
    & \quad + \bigg( \sum_{k=-\infty}^{{M-2}} \bigg(\int_{x_{k+1}}^b w(t) \bigg(\int_t^{b} u \bigg)^{\frac{r}{q}} \,dt  \bigg)^{\frac{1}{1-r}} \int_{x_k}^{x_{k+1}}  d\big[V_p(a,t)^{\frac{r}{1-r}}\big] \bigg)^{\frac{1-r}{r}} \\
& \lesssim  B_2^* + \bigg( \sum_{k=-\infty}^{{M-1}} \int_{x_k}^{x_{k+1}} \bigg(\int_t^{x_{k+1}} w(s) \bigg(\int_s^{b} u \bigg)^{\frac{r}{q}} \,ds  \bigg)^{\frac{1}{1-r}} d\big[V_p(a,t)^{\frac{r}{1-r}}\big] \bigg)^{\frac{1-r}{r}}.
\end{align*}
By integrating by parts again, we obtain that
\begin{align}
C_5 &  \lesssim  B_2^* + \bigg( \sum_{k=-\infty}^{{M-1}} \int_{x_k}^{x_{k+1}} \bigg(\int_t^{x_{k+1}} w(s) \bigg(\int_s^{b} u \bigg)^{\frac{r}{q}} \,ds  \bigg)^{\frac{r}{1-r}} w(t) \bigg(\int_t^{b} u \bigg)^{\frac{r}{q}} V_p(a, t)^{\frac{r}{1-r}} \,dt \bigg)^{\frac{1-r}{r}}\notag\\
& \leq B_2^* + \bigg( \sum_{k=-\infty}^{{M-1}} \int_{x_k}^{x_{k+1}} \bigg(\int_t^{x_{k+1}} w  \bigg)^{\frac{r}{1-r}} w(t) \bigg(\int_t^{b} u \bigg)^{\frac{r}{q(1-r)}} V_p(a, t)^{\frac{r}{1-r}} \,dt \bigg)^{\frac{1-r}{r}}\notag\\
& \leq B_2^* + \bigg( \sum_{k=-\infty}^{{M-1}} \int_{x_k}^{x_{k+1}} \bigg(\int_t^{x_{k+1}} w  \bigg)^{\frac{r}{1-r}} w(t) \bigg( \esup_{s \in (t, b)} \bigg(\int_s^{b} u \bigg)^{\frac{r}{q(1-r)}} V_p(a, s)^{\frac{r}{1-r}} \bigg) \,dt \bigg)^{\frac{1-r}{r}}\notag\\
& \lesssim B_2^* + \tilde{A}_3. \label{C5<B2*+A3}
\end{align}
Consequently, we arrive at $\tilde{A}_3 + B_2^* \lesssim C_4 + C_5 \lesssim  \tilde{A}_3 + B_2^*$.

\rm(iv) Let $p \leq 1 $, $r < 1$, $q < 1$. We know from Theorem~\ref{T:disc-p<1}(iv) that $C \approx A_4^* + B_2^*$. First we will prove that $A_4^* + B_2^* \approx \tilde{A}_4 + B_2^*$, where
\begin{equation*}
\tilde{A}_4 :=\bigg( \sum_{k=-\infty}^{{M-1}} 2^{\frac{k}{1-r}} \bigg( \int_{x_k}^b \bigg(\int_{t}^{b} u \bigg)^{\frac{q}{1-q}} u(t) V_p(a, t)^{\frac{q}{1-q}} \,dt \bigg)^{\frac{r(1-q)}{q(1-r)}} \bigg)^{\frac{1-r}{r}}.
\end{equation*}

It is clear that $A_4^* \leq \tilde{A}_4$. On the other hand, since $\max\{A_4^*, B_2^*\}<\infty$ implies that $\max\{A_2^*, B_1^*\}<\infty$, by using the same argument we applied in case (ii), \eqref{IBP} holds. Then, \eqref{inc-sum-sum} combined with  \eqref{IBP} and \eqref{limit<supremum} yields that
\begin{align*}
\tilde{A}_4 & = \bigg( \sum_{k=-\infty}^{{M-1}} 2^{\frac{k}{1-r}} \bigg( \sum_{i=k}^{M-1}\int_{x_i}^{x_{i+1}} \bigg(\int_t^{b} u \bigg)^{\frac{q}{1-q}}u(t) V_p(a, t)^{\frac{q}{1-q}} \,dt \bigg)^{\frac{r(1-q)}{q(1-r)}} \bigg)^{\frac{1-r}{r}}\\
& \approx \bigg( \sum_{k=-\infty}^{{M-1}} 2^{\frac{k}{1-r}} \bigg( \int_{x_k}^{x_{k+1}} \bigg(\int_t^{b} u \bigg)^{\frac{q}{1-q}}u(t) V_p(a, t)^{\frac{q}{1-q}} \,dt \bigg)^{\frac{r(1-q)}{q(1-r)}} \bigg)^{\frac{1-r}{r}}\\
&\lesssim \bigg( \sum_{k=-\infty}^{{M-1}} 2^{\frac{k}{1-r}} \bigg( \int_{x_k}^{x_{k+1}} \bigg(\int_{t}^{b} u \bigg)^{\frac{1}{1-q}} d\big[V_p(a, t)^{\frac{q}{1-q}}\big] \bigg)^{\frac{r(1-q)}{q(1-r)}} \bigg)^{\frac{1-r}{r}}\\
    &\quad +\bigg( \sum_{k=-\infty}^{{M-1}} 2^{\frac{k}{1-r}} \esup_{t \in (x_k, x_{k+1})} \bigg(\int_{t}^{b} u \bigg)^{\frac{r}{q(1-r)}} V_p(a, t)^{\frac{r}{1-r}}\bigg)^{\frac{1-r}{r}}\\
& \lesssim \bigg( \sum_{k=-\infty}^{{M-1}} 2^{\frac{k}{1-r}} \bigg( \int_{x_k}^{x_{k+1}} \bigg(\int_{t}^{x_{k+1}} u \bigg)^{\frac{1}{1-q}} d\big[V_p(a, t)^{\frac{q}{1-q}}\big] \bigg)^{\frac{r(1-q)}{q(1-r)}} \bigg)^{\frac{1-r}{r}}\\
    &\quad +\bigg( \sum_{k=-\infty}^{{M-1}} 2^{\frac{k}{1-r}} \esup_{t \in (x_k, x_{k+1})} \bigg(\int_{t}^{x_{k+1}} u \bigg)^{\frac{r}{q(1-r)}} V_p(a, t)^{\frac{r}{1-r}}\bigg)^{\frac{1-r}{r}}\\
    & \quad + \bigg( \sum_{k=-\infty}^{{M-2}} 2^{\frac{k}{1-r}} \bigg(\int_{x_{k+1}}^b u \bigg)^{\frac{r}{q(1-r)}} V_p(a, {x_{k+1}})^{\frac{r}{1-r}}\bigg)^{\frac{1-r}{r}}.
\end{align*}
Integration by parts gives
\begin{align*}
\tilde{A}_4 & \lesssim \bigg( \sum_{k=-\infty}^{{M-1}} 2^{\frac{k}{1-r}} \bigg( \int_{x_k}^{x_{k+1}} \bigg(\int_{t}^{x_{k+1}} u \bigg)^{\frac{q}{1-q}} u(t) V_p(a, t)^{\frac{q}{1-q}} \,dt \bigg)^{\frac{r(1-q)}{q(1-r)}} \bigg)^{\frac{1-r}{r}} \\
    &\quad +\bigg( \sum_{k=-\infty}^{{M-1}} 2^{\frac{k}{1-r}} \esup_{t \in (x_k, x_{k+1})} \bigg(\int_{t}^{x_{k+1}} u \bigg)^{\frac{r}{q(1-r)}} V_p(a, t)^{\frac{r}{1-r}}\bigg)^{\frac{1-r}{r}}\\
    & \quad + \bigg( \sum_{k=-\infty}^{{M-2}} 2^{\frac{k}{1-r}} \bigg(\int_{x_{k+1}}^b u \bigg)^{\frac{r}{q(1-r)}} V_p(a, {x_{k+1}})^{\frac{r}{1-r}}\bigg)^{\frac{1-r}{r}}\\
& =: \tilde{A}_{4,1} +\tilde{A}_{4,2}+\tilde{A}_{4,3}.
\end{align*}
It is easy to see that $\tilde{A}_{4,3} \lesssim B_2^*$. On the other hand,  observe that \eqref{upper-A1-A2} yields $\tilde{A}_{4,2} \lesssim \tilde{A}_{4,1}$. Moreover, using \eqref{V-cut}, we have
\begin{align*}
\tilde{A}_{4,1} & \lesssim A_4^* +  \bigg( \sum_{k=-\infty}^{{M-1}} 2^{\frac{k}{1-r}} \bigg( \int_{x_k}^{x_{k+1}} u \bigg)^{\frac{r}{q(1-r)}} V_p(a, x_k)^{\frac{r}{1-r}} \bigg)^{\frac{1-r}{r}} \\
& \lesssim A_4^* + B_2^*.
\end{align*}
Thus, we arrive at $\tilde{A}_4 \lesssim A_4^* + B_2^*$.

We proceed by proving $\tilde{A}_4 + B_2^* \approx C_5 + C_4$. It is clear by using \eqref{int.equiv} with $\alpha=\frac{r}{1-r}$ that $\tilde{A}_4 \approx C_6$. On the other hand, using \eqref{upper-A1-A2}, we conclude that $C_4 \lesssim C_6$ and $\tilde{A}_3 \lesssim \tilde{A}_4$. Thus, taking \eqref{B2*<C5} into consideration, we have $\tilde{A}_3 + B_2^* \lesssim \tilde{A}_4 + B_2^* \lesssim C_5 + C_4 \lesssim C_5 + C_6$. It remains to prove that $C_5 \lesssim \tilde{A}_4 + B_2^* $. We have already proved in \eqref{C5<B2*+A3} that $C_5 \lesssim \tilde{A}_3 + B_2^*$. Combining these estimates we arrive at $\tilde{A}_4 + B_2^* \lesssim C_5 + C_6 \lesssim \tilde{A}_4 + B_2^*$, and the proof is complete.
\end{proof}



\end{document}